\documentclass{amsart}
\usepackage{amsmath,amsthm,amssymb}

\usepackage{graphicx}

\theoremstyle{plain}

\newtheorem{theorem}{Theorem}[section]
\newtheorem{lemma}{Lemma}[section]

\newtheorem{condition}{Condition}[section]

\theoremstyle{definition}

\newtheorem{definition}{Definition}[section]
\newtheorem{example}{Example}[section]
\newtheorem{remark}{Remark}[section]

\numberwithin{equation}{section}

\newcommand{\ophi}{{\overline\varphi}}

\renewcommand{\phi}{{\varphi}}

\renewcommand{\phi}{{\varphi}}

\newcommand{\cV}{{\mathcal V}}

\newcommand{\cR}{{\mathcal R}}
\newcommand{\cH}{{\mathcal H}}

\newcommand{\cW}{{\mathcal W}}

\newcommand{\cU}{{\mathcal U}}

\newcommand{\oQ}{{\overline Q}}

\newcommand{\ob}{{\overline b}}

\newcommand{\hu}{{\hat u}}
\newcommand{\hv}{{\hat v}}

\newcommand{\bbR}{{\mathbb R}}
\newcommand{\bbN}{{\mathbb N}}

\renewcommand{\phi}{{\varphi}}

\begin{document}

\title{Systems of reaction-diffusion equations with spatially
distributed hysteresis}

\author{Pavel  Gurevich, \\
        Sergey  Tikhomirov}


\begin{abstract}
We study systems of reaction-diffusion equations with
discontinuous spatially distributed hysteresis in the right-hand
side. The input of hysteresis is given by a vector-valued function
of space and time. Such systems describe hysteretic interaction of
non-diffusive (bacteria, cells, etc.) and diffusive (nutrient,
proteins, etc.) substances leading to formation of spatial
patterns. We provide sufficient conditions under which the problem
is well posed in spite of the discontinuity of hysteresis. These
conditions are formulated in terms of geometry of manifolds
defining  hysteresis thresholds and the graph of initial data.
\end{abstract}



\thanks{
The research of the first author was supported by the DFG project
SFB 910. The research of the second author was supported by
Alexander von Humboldt Foundation, by the Chebyshev Laboratory
(Department of Mathematics and Mechanics, St. Petersburg State
University) under RF Government grant 11.G34.31.0026 and JSC ``Gazprom Neft'' and by posdoctoral fellowship of the Max Planck Institute for
Mathematics in the Sciences.}


\maketitle

\section{Setting of the problem}

\subsection{Introduction and setting} Reaction-diffusion equations
with spatially distributed hysteresis were first introduced
in~\cite{Jaeger1} to describe the growth of a colony of bacteria
(Salmonella typhimurium) and explain emerging spatial patterns of
the bacteria density. In~\cite{Jaeger1, Jaeger2}, numerical
analysis of the problem was carried out, however without rigorous
justification. First analytical results were obtained
in~\cite{Alt, VisintinSpatHyst} (see
also~\cite{AikiKopfova,Kopfova}), where existence of solutions for
multi-valued hysteresis was proved. Formal asymptotic expansions
of solutions were recently obtained in a special case
in~\cite{Ilin}. Questions about the uniqueness of solutions and
their continuous dependence on initial data as well as a thorough
analysis of pattern formation remained open. In this paper, we
formulate sufficient conditions that guarantee existence,
uniqueness, and continuous dependence of solutions on initial data
for systems of reaction-diffusion equations with discontinuous
spatially distributed hysteresis. Analogous conditions for scalar
equations have been considered by the authors
in~\cite{GurTikhSIAM13, GurTikhNonlinAnal}.

Denote $Q_T=(0,1)\times(0,T)$, where $T>0$. Let $\mathcal U
\subset\mathbb R^k$ and $\mathcal V \subset\mathbb R^l$
($k,l\in\mathbb N$) be closed sets. We assume throughout that
$(x,t)\in \overline Q_T$, $u(x,t)\in\mathcal U$,
$v(x,t)\in\mathcal V$.

We consider the system of reaction-diffusion equations
\begin{equation}\label{eq1}\left\{
\begin{aligned}
&u_t=D u_{xx}+f(u,v,W(\xi_0,u)),\\
&v_t=g(u,v,W(\xi_0,u))\\
\end{aligned}\right.
\end{equation}
with the initial and boundary conditions
\begin{equation}
u|_{t=0}=\phi(x),\quad v|_{t=0}=\psi(x),\quad
u_x|_{x=0}=u_x|_{x=1}=0. \label{eq3}
\end{equation}
Here $D$ is a positive-definite diagonal matrix; $W$ is a
hysteresis operator which maps an initial configuration function
$\xi_0(x)$ ($\in\{1,-1\}$) and an input function $u(x,\cdot)$
 to an output function $W(\xi_0(x),u(x,\cdot))(t)$. As
a function of $(x,t)$, $W(\xi_0,u)$ takes values in a set
$\mathcal W\subset \mathbb R^m$ ($m\in\mathbb N$). Now we shall
define this operator in detail.

Let $\Gamma_\alpha,\Gamma_\beta\subset\mathcal U$ be two disjoint
smooth manifolds of codimension one without boundary (hysteresis
``thresholds''). For simplicity, we assume that they are given by
$\gamma_\alpha(u)=0$ and $\gamma_\beta(u)=0$ with $\nabla
\gamma_\alpha(u)\ne 0$ and $\nabla \gamma_\beta(u)\ne 0$,
respectively, where $\gamma_\alpha$ and $\gamma_\beta$ are
$C^\infty$-smooth functions (in the general situation, atlases can
be used).

\begin{figure}[ht]
        \centering
      \includegraphics[width=6cm]{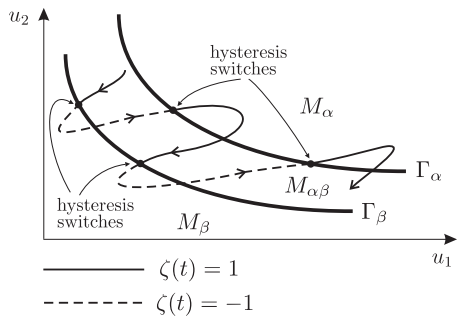}
        \caption{Regions of different behavior of hysteresis $W$}
        \label{fig1}
\end{figure}

Denote $M_\alpha=\{u\in\mathcal U:\gamma_\alpha(u)\le 0\}$,
$M_\beta=\{u\in\mathcal U:\gamma_\beta(u)\le 0\}$,
$M_{\alpha\beta}=\{u\in\mathcal U:\gamma_\alpha(u)>0,\
\gamma_\beta(u)>0\}$. Assume that
$M_\alpha\cap\Gamma_\beta=\varnothing$ and
$M_\beta\cap\Gamma_\alpha=\varnothing$ (Fig.~\ref{fig1}).

Next, we introduce locally H\"older continuous functions ({\it
hysteresis ``branches''})
$$
W_1: D(W_1)=M_\alpha\cup\overline M_{\alpha\beta}\to\mathcal
W,\qquad W_{-1}: D(W_{-1})=M_\beta\cup\overline M_{\alpha\beta}
\to\mathcal W.
$$

We fix $T>0$ and denote by $C_r[0,T )$ the   space of functions
continuous on the right in $[0,T )$. For any $\zeta_0\in\{1,-1\}$
({\it initial configuration}) and $u_0\in C([0,T];\mathcal U)$
({\it input}), we introduce the {\it configuration} function
$$
\zeta:\{1,-1\}\times C([0,T];\mathcal U)\to C_r[0,T),\quad
\zeta(t)=\zeta(\zeta_0,u_0)(t)
$$
 as follows. Let
$X_t=\{s\in(0,t]:u_0(s)\in\Gamma_\alpha\cup \Gamma_\beta\}$. Then
$ \zeta(0)=1$  if $u_0(0 )\in M_\alpha$, $\zeta(0)=-1$ if
$u_0(0)\in M_\beta$, $\zeta(0)=\zeta_0$ if $u_0(0)\in
M_{\alpha\beta}$;  for $t\in(0,T]$, $\zeta(t)=\zeta(0)$ if
$X_t=\varnothing$, $\zeta(t)= 1$ if $X_t\ne\varnothing$ and
$u_0(\max X_t)\in\Gamma_\alpha$, $\zeta(t)=-1$  if
$X_t\ne\varnothing$ and $u_0(\max X_t)\in\Gamma_\beta$
(Fig.~\ref{fig1}).

Now we introduce the {\it hysteresis operator\/}
 $ W: \{1,-1\}\times
C([0,T];\mathcal U)\to C_r[0,T ) $ by the following rule
(cf.~\cite{KrasnBook,Visintin,Klein}). For any initial
configuration $\zeta_0\in\{1,-1\}$  and input $u_0\in
C([0,T];\mathcal U)$, the function
$W(\zeta_0,u_0):[0,T]\to\mathcal W$ ({\it output}) is given by
\begin{equation}\label{eqHystTime}
W(\zeta_0,u_0)(t)=W_{\zeta(t)}(u_0(t)),
\end{equation}
where $\zeta(t)$ is the configuration function defined above.

Assume that the initial configuration and the input function
depend on spatial variable $x$. Denote them by $\xi_0(x)$ and
$u(x,t)$, where $u(x,\cdot)\in C([0,T];\mathcal U)$.
Using~\eqref{eqHystTime} and treating $x$   as a parameter, we
define the {\it spatially distributed hysteresis}
\begin{equation}\label{eqSDH}
W(\xi_0(x),u(x,\cdot))(t)=W_{\xi(x,t)}(u(x,t)),
\end{equation}
where    $ \xi(x,t)=\zeta(\xi_0(x),u(x,\cdot))(t) $ is the {\it
spatial configuration}.

\subsection{Functional spaces.} Denote by
$L_q(0,1)$, $q>1$, the standard Lebesgue space and by $W_q^l(0,1)$
with natural $l$ the standard Sobolev space. For a noninteger
$l>0$, denote by $W_q^l(0,1)$ the Sobolev space with the norm
$$
\|\phi\|_{W_q^{[l]}(0,1)}+\left(\int_0^1dx\int_0^1
\dfrac{|\phi^{([l])}(x)-\phi^{([l])}(y)|^q}{|x-y|^{1+q(l-[l])}}dy\right)^{1/q},
$$
where $[l]$ is the integer part of $l$. Introduce
 the anisotropic Sobolev spaces $ W_q^{2,1}(Q_T) $ with the norm
$ \left(\int_0^T \|u(\cdot,t)\|_{W_q^2(0,1)}^q\, dt+ \int_0^T
\|u_t(\cdot,t)\|_{L_q(0,1)}^q\, dt\right)^{1/q} $ and the space
$W_\infty^{0,1}(Q_T)$ of $L_\infty(0,1)$-valued functions
continuously differentiable on $[0,T]$  with the norm $
\|v\|_{L_\infty(Q_T)}+\|v_t\|_{L_\infty(Q_T)}.$ Denote by
$C^\gamma(\oQ_T)$, $\gamma\in(0,1)$, the H\"older space.

For the vector-valued functions, we  use the following notation.
If, e.g., $u(x,t)\in\mathcal U$ and each component of $u$ belongs
to $W_q^{2,1}(Q_T)$, then we  write  $u\in W_q^{2,1}(Q_T;\mathcal
U)$.

Throughout, we fix $q$ and $\gamma$ such that $ q\in (3,\infty)$
and $\gamma\in (0,1-3/q).$ This implies that $u,u_x\in
C^\gamma(\overline Q_T;\mathcal U)$ for $u\in
W_q^{2,1}(Q_T;\mathcal U)$ (see Lemma 3.3 in~\cite[Chap.
2]{LadSolUral}).

To define the space of initial data, we  use the fact that if
$u\in W_q^{2,1}(Q_T;\mathcal U)$, then the trace $u|_{t=0}$ is
well defined and belongs to $W_{q}^{2-2/q}((0,1);\mathcal U)$ (see
Lemma 2.4 in~\cite[Chap. 2]{LadSolUral}). We denote the norm in
the latter space by $\|\cdot\|_q$. Moreover, one can define the
space $W_{q,N}^{2-2/q}((0,1);\mathcal U)$ as the subspace of
functions from $W_q^{2-2/q}((0,1);\mathcal U)$ with the zero
Neumann boundary conditions.

We assume  that $\phi\in W_{q,N}^{2-2/q}((0,1);\mathcal U)$ and
$\psi\in L_\infty((0,1);\mathcal V)$ in~\eqref{eq3}.

\begin{definition}\label{defSolution}
A pair $(u,v)\in W_q^{2,1}(Q_T;\mathcal U)\times
W_\infty^{0,1}(Q_T;\mathcal V)$ is  a {\em solution of
problem~\eqref{eq1}, \eqref{eq3}}  if $W(\xi_0,u)$ is measurable
with respect to $(x,t)$ and~\eqref{eq1}, \eqref{eq3} hold.
\end{definition}

\subsection{Spatial transversality}
We will deal with the case where $\xi_0(x)$ has one discontinuity
point. Generalization to  finitely many discontinuity points is
straightforward.
\begin{condition}\label{condA1'}
\begin{enumerate}
\item For some $\ob\in(0,1)$, we have
\begin{equation}\label{eqxi0}
\xi_0(x)=1\ (x\le\ob),\quad \xi_0(x)=-1\ (x>\ob).
\end{equation}
\item For $x \in [0,\ob]$, we have $\phi(x) \in M_{\alpha} \cup
M_{\alpha \beta}$ or, equivalently, $\gamma_\beta(\phi(x))>0$.

\item For $x \in (\ob,1]$, we have $\phi(x) \in M_{\beta} \cup
M_{\alpha \beta}$ or, equivalently, $\gamma_\alpha(\phi(x)) > 0$.

\item If $\gamma_\alpha(\phi(\ob))=0$, then
$\dfrac{d}{dx}\gamma_\alpha(\phi(x))\Big|_{x=\overline b}>0$.

\end{enumerate}
\end{condition}
It follows from Condition~\ref{condA1'} that the hysteresis
in~\eqref{eqSDH} at the initial moment  equals $W_1(\phi(x))$ for
$x\le\ob$ and $W_2(\phi(x))$ for $x>\ob$. Items 2 and 3 in
Condition~\ref{condA1'} are necessary for the hysteresis to be
well-defined at the initial moment, while item 4 is an essential
assumption. We   refer to item 4 as the {\it spatial
transversality} and  say that $\phi(x)$ {\it is transverse with
respect to the spatial configuration $\xi_0(x)$}. This means that
if   $\phi(\overline b)\in\Gamma_\alpha$, then the vector
$\phi'(\overline b)$ is transverse to the manifold $\Gamma_\alpha$
at this point.


%
%
%
%

Consider  time-dependent functions $u(x,t)$ such that $u,u_x\in
C(\overline Q_T;\mathcal U)$.

\begin{definition}
We say that a function $u$ is {\em transverse} {\em on} $[0,T]$
$($with respect to a spatial configuration $\xi(x,t))$ if, for
every fixed $t\in[0,T]$,  either $\xi(\cdot,t)$ has no
discontinuity points for $x\in(0,1)$, or it has one discontinuity
point and the function $u(\cdot,t)$ is transverse with respect to
the spatial configuration $\xi(\cdot,t)$.
\end{definition}

\begin{definition}
A function $u$ {\em preserves  spatial topology $($of a spatial
configuration $\xi(x,t))$ on $[0,T]$} if there is $M>0$ such that,
for $t\in[0,T]$, there is a continuous function $b(t)\in (0,1)$
such that $\xi(x,t)=1$ for $x\le b(t)$ and $\xi(x,t)=-1$ for
$x>b(t)$.
\end{definition}

The solution from Definition~\ref{defSolution} is called {\em
transverse} ({\em preserving spatial topology}) if the function
$u(x,t)$ is transverse (preserving spatial topology).

\begin{remark}
The function  $b(t)$ defining discontinuity of $\xi(x,t)$ plays a
role of free boundary, which has much in common with   free
boundary in parabolic obstacle problems (see, e.g., \cite{UrAp,
ShUr} and the references therein). However, in our case, the
behavior of $b(t)$ is defined differently.
\end{remark}

\subsection{Assumptions on the right-hand side.} First, we assume
the following.

\begin{condition}\label{condRegularity}
The functions $f(u,v,w)$ and $g(u,v,w)$ are locally Lipschitz
continuous in $\mathbb R^k\times\mathbb R^l\times\mathbb R^m$.
\end{condition}

Next, we formulate  dissipativity conditions for $f$ and $g$.

In the following condition, we denote by $\mathcal U_\mu$, $\mu\ge
0$, closed parallelepipeds in $\mathcal U$ with the edges parallel
to respective coordinate axes such that $\phi(x)\in\mathcal U_\mu$
for all $x\in[0,1]$.


\begin{condition}\label{condDissipativity-f}
 There is a parallelepiped $\mathcal U_0$ and, for each
sufficiently small $\mu>0$, there is a parallelepiped $\mathcal
U_\mu$ and a locally Lipschitz continuous function $f_\mu(u,v)$
such that
\begin{enumerate}
\item $|f_\mu(u,v)|$ converges to $0$ uniformly on compact sets in
$\mathcal U \times \mathcal V$ as $\mu\to 0$,

\item At each point $u\in\partial \cU_0\cap D(W_{\pm1})$,
$v\in\mathcal V$,  the vector $f(u,v,W_{\pm1}(u))+f_\mu(u,v)$
points strictly inside $\mathcal U_0$.

\item At each point $u\in\partial \cU_\mu\cap D(W_{\pm1})$,
$v\in\mathcal V$,        the vector
$f(u,v,W_{\pm1}(u_\mu))+f_\mu(u,v)$ points strictly inside
    $\mathcal U_\mu$ for all $u_\mu\in \mathcal U_\mu$.
\end{enumerate}
\end{condition}

To formulate the  assumption on $g$, we fix $\cU_0$ satisfying
Condition~\ref{condDissipativity-f} and set
\begin{equation}\label{eqcW0}
\cW_0=\bigcup_{j=\pm 1}\{W_j(u):u\in\cU_0\}.
\end{equation}

\begin{condition}\label{condDissipativity-g}
For any $T_0>0$, there is a compact $\mathcal V_0=\mathcal
V_0(T_0,\cU_0)\subset\mathcal V$ such that $\psi(x)\in\cV_0$
$((x\in(0,1))$ and  the Cauchy problem
\begin{equation}\label{eqODEv}
 v_t=g(u_0(x,t),v,w_0(x,t)),\quad v|_{t=0}=\psi(x)
\end{equation}
has a solution $v\in W_\infty^{0,1}(Q_{T_0};\bbR^l)$ satisfying
$v(x,t)\in\cV_0$  whenever
$$
\begin{gathered}
u_0\in L_\infty(Q_{T_0};\cU),\quad w_0\in
L_\infty(Q_{T_0};\cW),\\
u_0(x,t)\in\cU_0, \quad w_0(x,t)\in\cW_0\quad ((x,t)\in Q_{T_0}).
\end{gathered}
$$
\end{condition}

\begin{remark}
It follows from~\cite[Theorem 1, p. 111]{Rothe} that
system~\eqref{eqODEv}
 has a unique  solution $v\in W_\infty^{0,1}(Q_{T_0};\bbR^l)$ for a sufficiently small
 $T_0>0$.  Condition~\ref{condDissipativity-g} additionally guarantees
the absence of blow-up.

In particular, the uniform boundedness of $v$ holds if
$|g(u,v,w)|\le A(u,w)|v|+B(u,w)$, where $A(u,w)$ and $B(u,w)$ are
bounded on compact sets (see Example~\ref{exBacteria}). However,
if $\mathcal V\ne\mathbb R^l$, one must additionally check that
$v$ never  leaves $\mathcal V$. To fulfill
Condition~\ref{condDissipativity-g}, one could alternatively
assume the existence of invariant parallelepiped for $g$
(similarly to Condition~\ref{condDissipativity-f}).
\end{remark}

\begin{example}\label{exBacteria}
The hysteresis operator and the right-hand side in the present
paper apply to a model describing a growth of a colony of bacteria
(Salmonella typhimurium) on a petri plate (see~\cite{Jaeger1,
Jaeger2}). Let $u_1(x,t)$ and $u_2(x,t)$ denote the concentrations
of diffusing buffer (pH level) and histidine (nutrient),
respectively, while $v(x,t)$ denote the density of nondiffusing
bacteria. These three unknown functions satisfy the following
equations:
\begin{equation}\label{eqBact1}
\left\{
\begin{aligned}
u_{1t}&=D_1\Delta u_1-a_1 W(\xi_0,u)v,\\
u_{2t}&=D_2\Delta u_2-a_2 W(\xi_0,u)v,\\
v_t&=aW(\xi_0,u)v,\\
\end{aligned}\right.
\end{equation}
supplemented by initial and no-flux (Neumann) boundary conditions.
In~\eqref{eqBact1}, $D_1,D_2,a,a_1,a_2>0$ are given constants and
$W(\xi_0,u)$ is the hysteresis operator. In this example, we have
$\mathcal U=\{u\in\mathbb R^2:u_1,u_2\ge 0\}$, $\mathcal
V=[0,\infty)$, $\mathcal W=[0,\infty)$. The hysteresis thresholds
$\Gamma_\alpha$ and $\Gamma_\beta$ are the curves on the plane
given by $\gamma_\alpha(u):=-u_1+a_\alpha/u_2+b_\alpha=0$ and
$\gamma_\beta(u):=u_1-a_\beta/u_2-b_\beta=0$, respectively, where
$a_\alpha,a_\beta,b_\alpha,b_\beta>0$ are some constants
(Fig.~\ref{fig1}); the hysteresis ``branches'' are given by
functions $W_1(u)$ $(>0)$ and $W_{-1}(u)$ $(\equiv 0)$.
\end{example}

\section{Main results.} In what follows, we assume that
Conditions~$\ref{condA1'}$--$\ref{condDissipativity-g}$ hold.

\begin{theorem}[local existence]\label{thLocalExistence}
There is a number $T>0$ such that
\begin{enumerate}
\item There is at least one solution of problem~\eqref{eq1},
\eqref{eq3} in~$Q_T;$

\item Any solution in $Q_T$ is transverse and preserves spatial
topology.
\end{enumerate}
\end{theorem}

\begin{theorem}[continuation]\label{tCont}
Let $(u,v)$ be a transverse topology preserving
 solution of problem~\eqref{eq1}, \eqref{eq3}  in $Q_T$ for some $T>0$. Then
it can be continued to an   interval $[0,T_{max})$, where
$T_{max}>T$ has the following properties. 1. For any
$t_0<T_{max}$, the pair $(u, v)$ is a transverse solution of
problem~\eqref{eq1}, \eqref{eq3}  in $Q_{t_0}$. 2. Either
$T_{max}=\infty$, or  $T_{max}<\infty$ and $(u,v)$ is a solution
in $Q_{T_{max}}$, but $u(\cdot, T_{max})$ is not transverse with
respect to $\xi(\cdot,T_{max})$.
\end{theorem}

\begin{theorem}[continuous dependence on initial data]\label{tContDepInitData}
Assume the following.
\begin{enumerate}
\item There is a number $T>0$ such that problem~\eqref{eq1},
\eqref{eq3} with initial functions $\phi,\psi$ and initial
configuration $\xi_0(x)$ defined by its discontinuity point $\ob$
admits a unique transverse topology preserving solution $(u,v)$ in
$Q_s$ for any $s\le T$.

\item Let $\phi_n\in W_{q,N}^{2-2/q}((0,1);\mathcal U)$,
$\psi_n\in L_\infty((0,1);\mathcal V)$, $n=1,2,\dots$, be a
sequence of other initial functions such that
$\|\phi-\phi_n\|_q\to 0$, $\|\psi-\psi_n\|_{L_q((0,1);\mathcal
V)}\to 0$ as $n\to\infty$.

\item Let $\xi_{0n}(x)$,  $n=1,2,\dots$, be a sequence of other
initial configurations defined by their discontinuity points
$\ob_{n}$ similarly to~\eqref{eqxi0}  and $\ob_{n}\to\ob$ as
$n\to\infty$.
\end{enumerate}
Then, for all sufficiently large $n$, problem~\eqref{eq1},
\eqref{eq3} with initial data $(\phi_n,\psi_n,\xi_{0n})$ has at
least one transverse topology preserving solution $(u_n,v_n)$.
Each sequence of such solutions satisfies
$$
\begin{gathered}
\|u_n-u\|_{W_q^{2,1}(Q_T;\mathcal U)}\to 0,\quad
 \|b_{jn}-b_j\|_{C[0,T]}\to 0,\\
\sup\limits_{t\in[0,T]}\left(\|v_n(\cdot,t)-
v(\cdot,t)\|_{L_q((0,1);\cV)}+\|v_{nt}(\cdot,t)-
v_t(\cdot,t)\|_{L_q((0,1);\cV)}\right)\to 0
\end{gathered}
$$
as $n\to\infty,$ where $b(t)$ and $b_{n}(t)$ are the respective
discontinuity points of the configuration functions $\xi(x,t)$ and
$\xi_n(x,t)$.
\end{theorem}

\begin{remark}
If one a priori knows that all $u_n$ are transverse on some
interval $[0,T]\subset[0,T_{max})$, then one can prove that $u_n$
approximate $u$ on $[0,T]$ even if $u$ is not topology preserving
on $[0,T]$.
\end{remark}



Now we discuss the uniqueness of solutions.   We strengthen the
assumption about local H\"older continuity of $W_{\pm1}$. Let
$\mathcal U_0$ be the set from
Condition~\ref{condDissipativity-f}.
\begin{condition}\label{condUniq}
There are numbers $K>0$ and $\sigma\in[0,1)$  such that
$$
\begin{aligned}
|W_1(u)-W_1(\hat u)|&\le
\dfrac{K}{(\gamma_\beta(u))^\sigma+(\gamma_\beta(\hat
u))^\sigma}|u-\hat u|\quad & &\forall u,\hat u\in M_\alpha\cup\overline M_{\alpha\beta},\\
|W_{-1}(u)-W_{-1}(\hat u)|&\le
\dfrac{K}{(\gamma_\alpha(u))^\sigma+(\gamma_\alpha(\hat
u))^\sigma}|u-\hat u| \quad & &\forall u,\hat u\in
M_\beta\cup\overline M_{\alpha\beta}.
\end{aligned}
$$
\end{condition}
We refer readers to \cite{GurTikhNonlinAnal} for the discussion about functions satisfying this condition.

\begin{theorem}[uniqueness]\label{tUniqueness}
Assume additionally that Condition~$\ref{condUniq}$ holds. Let
$(u,v)$ and $(\hat u,\hat v)$ be two transverse solutions of
problem~\eqref{eq1}, \eqref{eq3} in $Q_T$ for some $T>0$. Then
$(u,v)=(\hat u,\hat v)$.
\end{theorem}

\section{Local existence, continuation\\ and continuous dependence of solutions on initial data}

In this section, we prove
Theorems~\ref{thLocalExistence}--\ref{tContDepInitData}.
Throughout the section, we fix $\cU_0$ satisfying
Condition~\ref{condDissipativity-f} and $\cW_0$ given
by~\eqref{eqcW0}. Next, we fix some $T_0\in(0,1]$ and then $\cV_0$
satisfying Condition~\ref{condDissipativity-g}.

\subsection{Preliminaries}
The following result is straightforward.
\begin{lemma}\label{lab}
\begin{enumerate}
\item Let $\lambda\in[0,1)$,  $a\in C^\lambda[0,T]$, and
$b(t)=\max\limits_{s\in[0,t]} a(s)$. Then $b\in C^\lambda[0,T]$
and $ \|b\|_{C^\lambda[0,T]}\le\|a\|_{C^\lambda[0,T]}. $

\item If $a_j\in C[0,T]$ and $b_j(t)=\max\limits_{s\in[0,t]}
a_j(s)$, $j=1,2$, then $ \|b_1-b_2\|_{C[0,T]}\le
\|a_1-a_2\|_{C[0,T]}. $
\end{enumerate}
\end{lemma}

For some $T\le T_0$, $u_1\in L_\infty(Q_T;\cU)$, and $b_1\in
C[0,T]$ such that $u_1(x,t)\in\cU_0$ ($(x,t)\in Q_T$) and
$b_1(t)\in [0,1)$ ($t\in[0,T]$), we define the function $w_1(x,t)$
by
\begin{equation}\label{eqv0}
w_1(x,t)=\begin{cases} W_1(u_1(x,t)), & 0\le x\le  b_1(t),\\
W_{-1}(u_1(x,t)), & b_1(t)<x\le 1;
\end{cases}
\end{equation}
here we assume $W_{\pm1}(u_1)$ to be extended to  $\cU_0$ without
loss of regularity.

\begin{lemma}\label{lu0b0v0}
Let $u_1,b_1$ be functions with the above properties, and let
$\hat u_1,\hat b_1$ be functions with the same properties. Let
$w_1$ be defined by~\eqref{eqv0} and $\hat w_1$ by~\eqref{eqv0}
with $\hat u_1$ and $\hat b_1$ instead of $u_1$ and $b_1$,
respectively. Then, for any $p\in[1,\infty)$,
\begin{equation*}\label{equ0b_10v0}
\|w_1-\hat w_1\|_{L_p(Q_T;\cW)}\le c_0 \left(T^{1/p} \|u_1-\hat
u_1\|_{L_\infty(Q_T;\cU)}^{\sigma_0}+\|b_1-\hat
b_1\|_{L_1(0,T)}^{1/p}\right),
\end{equation*}
where $\sigma_0$ is a H\"older exponent for $W_{\pm 1}(u_1)$ and
$c_0>0$ depends on $\cU_0$ and $p$, but does not depend on
$u_1,b_1,T$.
\end{lemma}
\proof We fix some $t\in[0,T]$ and assume that $b_1(t)\le \hat
b_1(t)$ for this $t$. Then, using~\eqref{eqv0} and omitting the
arguments of the integrands, we have
$$
\begin{aligned}
\int\limits_0^1 |w_1-\hat w_1|^p\,
dx&=\int\limits_0^{b_1(t)}|W_1(u_1)-W_1(\hat
u_1)|^p\,dx+\int\limits_{\hat b_1(t)}^1|W_{-1}(u_1)-W_{-1}(\hat
u_1)|^p\,dx\\
& +\int\limits_{b_1(t)}^{\hat b_1(t)}|W_{-1}(u_1)-W_1(\hat
u_1)|^p\,dx.
\end{aligned}
$$
Using the H\"older  continuity and the boundedness of $W_{\pm
1}(u_1)$ for $u_1\in\cU_0$ and  integrating with respect to $t$
from $0$ to $T$, we complete the proof.
\endproof

Now we introduce sets that ``measure'' the spatial transversality.
Denote by $E_m$, $m\in\mathbb N$, the set of triples
$(\phi,\psi,\xi_0)$ such that $\phi\in
W_{q,N}^{2-2/q}((0,1);\cU)$, $\psi\in L_\infty((0,1);\cV)$,
$\xi_0(x)$ is of the form~\eqref{eqxi0}, and the following hold:
\begin{enumerate}
\item $\ob \in [1/m, 1 - 1/m]$,

\item $\gamma_\beta(\phi(x))\ge 1/m^2$ for $x \in [0, \ob]$,

\item $\gamma_\alpha(\phi(x)) \ge 1/m^2$ for $x  \in [\ob + 1/m,
1]$,

\item if $x\in [\ob, \ob + 1/m]$ and $\gamma_\alpha(\phi(x)) \in
[0,  1/m^2]$, then $\dfrac{d}{dx}\gamma_\alpha(\phi(x))\ge 1/m$,

\item $\|\phi\|_q\le m$ and $\|\psi\|_{L_\infty((0,1);\cV)}\le m$.
\end{enumerate}

It is easy to check that  $E_m\subset E_{m+1}$. Moreover, one can
show (Lemma~2.25 in~\cite{GurTikhSIAM13}) that the union of all
sets $E_m$ coincides with the set of all  data satisfying
Condition~$\ref{condA1'}$. From now on, we fix $m\in\bbN$ such
that $(\phi,\psi,\xi_0)\in E_m$.

The next lemma follows from the implicit function theorem and
Lemma~\ref{lab}.

\begin{lemma}\label{lOneRoot}
Let $\lambda\in(0,1)$,   $u_1,u_{1x}\in C^\lambda(\oQ_{T_0};\cU)$,
$$
\|u_1\|_{C^\lambda(\oQ_{T_0};\cU)}+\|u_{1x}\|_{C^\lambda(\oQ_{T_0};\cU)}\le
c
$$
for some $c>0$, $u_1|_{t=0}=\varphi(x)$, and
$(\varphi,\psi,\xi_0)\in E_m$. Then there is
$T_1=T_1(m,\lambda,c)\le T_0$ and a natural number
$N_1=N_1(m,\lambda,c)\ge m$ which do not depend on $u,\phi,\xi_0$
such that the following is true for any $t\in[0,T_1]$.
\begin{enumerate}
\item The equation $\gamma_{\alpha}(u_1(x,t))=0$ for $x\in
[\ob,1]$ has no more than one root. If this root exists, we denote
it by $a_1(t);$ otherwise, we set $a_1(t)=\ob$. One has
$a_1(t)\in[\ob,\ob+1/N_1]$, $a_1\in C^\lambda[0,T_1]$.

\item The hysteresis $\cH(\xi_0,u_1)$  and its configuration
function $\xi_1(x,t)$ have  exactly one discontinuity point
$b_1(t);$ moreover, $b_1(t)=\max\limits_{s\in[0,t]}a_1(s)$,
$b_1\in C^\lambda[0,T_1]$.
\end{enumerate}
\end{lemma}

\subsection{Auxiliary problem}

Consider functions $u_1\in L_\infty(Q_T;\cU)$ and $w_1\in
L_\infty(Q_T;\cW)$ such that
$$
u_1(x,t)\in\cU_0,\quad w_1(x,t)\in\cW_0\quad ((x,t)\in Q_T)
$$
for some $T>0$. Define the functions
\begin{equation}\label{eqf0}
f_1(u,v,x,t)=f(u,v,w_1(x,t)),\quad
g_1(v,x,t)=g(u_1(x,t),v,w_1(x,t)).
\end{equation}
 Consider the auxiliary problem
\begin{equation}\label{eqPf0}\left\{
\begin{aligned}
&u_t=D u_{xx}+f_1(u,v,x,t),\\
&v_t=g_1(v,x,t)\\
&u|_{t=0}=\phi(x),\quad v|_{t=0}=\psi(x), \quad
u_x|_{x=0}=u_x|_{x=1}=0.
\end{aligned}\right.
\end{equation}
Set $f_U=\sup f(u,v,w)$ and $g_U=\sup g(u,v,w)$, where $(u,v,w)\in
\cU_0\times\cV_0\times\cW_0$.

The next result follows from the standard estimates for solutions
of linear parabolic equations~\cite{LadSolUral}, from
Conditions~\ref{condRegularity}--\ref{condDissipativity-g}
combined with the principle of invariant
rectangles~\cite{Smoller}, and from Lemma~\ref{lOneRoot}.

\begin{lemma}\label{ltNP0}
\begin{enumerate}
\item For any $T\le T_0$, problem~\eqref{eqPf0} has a unique
solution $(u,v)\in W_q^{2,1}(Q_T;\cU)\times
W_\infty^{0,1}(Q_T;\cV)$ and
\begin{gather}
\notag
u(x,t)\in\cU_0,\quad v(x,t)\in\cV_0\quad ((x,t)\in Q_T),\\
\notag
\|u\|_{W_q^{2,1}(Q_T;\cU)}+\max_{t\in[0,T]}\|u(\cdot,t)\|_{q} \le
c_1(\|\phi\|_{q}+f_U),\\
\notag \|v\|_{W_\infty^{0,1}(Q_T;\cV)}\le
\|\psi\|_{L_\infty((0,1);\cV)}+2g_U,\\
\label{eqNP02*}
 \|u\|_{C^\gamma(\oQ_T;\cU)}+\|u_x\|_{C^\gamma(\oQ_T;\cU)} \le
c_2 (\|\phi\|_{q}+f_U),
\end{gather}
where $c_1,c_2>0$ depend only on $T_0$.

\item If $u_n,v_n$, $n=1,2,\dots$, are solutions of
problem~\eqref{eqf0}, \eqref{eqPf0} with $u_1,w_1$ replaced by
$u_{1n},w_{1n}$ $($with the same properties$)$ and
\begin{equation*}
\|u_{1n}-u_1\|_{L_\infty(Q_T;\cU)}+\|w_{1n}-w_1\|_{L_q(Q_T;\cW)}\to
0,\quad n\to\infty,
\end{equation*}
 then
\begin{equation*}
\|u_n-u\|_{W_q^{2,1}(Q_T;\cU)}+\|v_n-v\|_{W_q^{0,1}(Q_T;\cV)}\to
0,\quad n\to\infty.
\end{equation*}

\item There is  $T_2=T_2(m)\le T_0$ and a natural number
$N_2=N_2(m)\ge m$ such that, for any $t\in[0,T_2]$, conclusions
(1) and (2) from Lemma~$\ref{lOneRoot}$ hold for $u(x,t)$, for the
corresponding ``root'' function $a(t)$, for the configuration
function $\xi(x,t)$ of the hysteresis $\cH(\xi_0,u)$, for its
discontinuity point $b(t)$, and for $T_2,N_2$ instead of
$T_1,N_1$. Furthermore, $(u(\cdot,t),v(\cdot,t), \xi(\cdot,t))\in
E_{N_2}$.
\end{enumerate}
 \end{lemma}

\subsection{Local existence: proof of Theorem~$\ref{thLocalExistence}$}

1. Let us prove the first assertion.

1.1.  Fix $\lambda$ in Lemma~\ref{lOneRoot} such that
$\lambda\in(0,\gamma)$. Fix $c_2$ from Lemma~\ref{ltNP0}. Set
$c=c_{\lambda,\gamma}c_2(m+f_U)$, where $c_{\lambda,\gamma}>0$ is
the embedding constant such that $\|u\|_{C^\lambda(\oQ_T;\cU)}\le
c_{\lambda,\gamma} \|u\|_{C^\gamma(\oQ_T;\cU)}$. Set
$T=\min(T_1,T_2)$, where $T_1, T_2$ are defined in
Lemmas~\ref{lOneRoot}, \ref{ltNP0}.

Let $R^\lambda(\oQ_T)$ be the set of functions $u(x,t)$ such that
$u|_{t=0}=\phi(x)$,
\begin{gather}
\notag u,u_x\in C^\lambda(\oQ_T;\cU),\quad u(x,t)\in\cU_0\quad
((x,t)\in
Q_T),\\
\label{eqProofLocalExistence1}
\|u\|_{C^\lambda(\oQ_T;\cU)}+\|u_x\|_{C^\lambda(\oQ_T;\cU)} \le c.
\end{gather}
 The set $R^\lambda(\oQ_T)$ is a closed convex subset of the
Banach space endowed with the norm given by the left-hand side
in~\eqref{eqProofLocalExistence1}. Similarly, we define
$R^\gamma(\oQ_T)$.

1.2. We construct a map $\cR:R^\lambda(\oQ_T)\to R^\gamma(\oQ_T)$.
Take any $u_1\in R^\lambda(\oQ_T)$ and define $a_1(t)$ and
$b_1(t)$ according to Lemma~\ref{lOneRoot}. Then define $w_1(x,t)$
by~\eqref{eqv0} and, using this $w_1$, define $f_1,g_1$
by~\eqref{eqf0}. Finally apply Lemma~\ref{ltNP0} and obtain a
solution $(u,v)$ of auxiliary problem~\eqref{eqPf0}. We now define
$\cR: u_1\mapsto u$.

The operator $\cR$ is continuous. Indeed, it is not difficult to
check that the mapping $u_1\mapsto a_1$ is continuous from
$R^\lambda(\oQ_T)$ to $C[0,T]$. Thus, the continuity of $\cR$
follows by consecutively applying Lemmas~\ref{lab} (part 2),
\ref{lu0b0v0}, \ref{ltNP0} (part 2), and the continuity of the
embedding $W_q^{2,1}(Q_T;\cV)\subset R^\gamma(\oQ_T)$.

Furthermore, due to~\eqref{eqNP02*} and the choice of $c$, the
operator $\cR$ maps $R^\lambda(\oQ_T)$ into itself. As an operator
acting from $R^\lambda(\oQ_T)$ into itself, it is compact due
to~\eqref{eqNP02*} and the compactness of the embedding
$R^\gamma(\oQ_T)\subset R^\lambda(\oQ_T)$. Therefore, applying the
Schauder fixed-point, we conclude the proof of the first assertion
of the theorem. Note that   $\cR$ is not Lipschitz continuous
(mind the exponent $1/p$ in~\eqref{equ0b_10v0}). Hence, the
contraction principle does not apply. We prove uniqueness
separately in Sec.~\ref{secUniqueness}.

2. The second assertion follows by applying the principle of
invariant rectangles (see~\cite{Smoller}) and Lemma~\ref{ltNP0}.


\subsection{Continuation: proof of Theorem~$\ref{tCont}$}
Theorem~\ref{tCont} follows from part~3 of Lemma~\ref{ltNP0} and
from the following fact (see Lemma~2.25 in~\cite{GurTikhSIAM13}).
Assume (1)   $(\phi_m,\psi_m,\xi_m)\in E_m\setminus E_{m-1}$,
$m=2,3,\dots$; (2) $\|\phi_m-\phi\|_q\to 0$ and
$\|\psi_m-\psi\|_{L_\infty((0,1);\cV)}\to 0$ as $m\to \infty$ for
some $\phi\in W_q^{2-2/q}((0,1);\mathcal U)$ and $\psi\in
L_\infty((0,1);\cV)$; (3) $ \ob_m-\ob\to0$ as $m\to \infty$ for
some $\ob\in[0,1]$. Then $\ob\in\{0,1\}$ or $\phi(x)$ is not
transverse with respect to $\xi_0(x)$, where $\xi_0(x)$ is given
by~\eqref{eqxi0}.

\subsection{Continuous dependence on initial data: proof of Theorem~$\ref{tContDepInitData}$}\label{subsubsecSmallInt1}
1. It suffices to prove the theorem for a sufficiently small time
interval. Since $(\phi,\psi,\xi_0)\in E_m$, it is easy to show
that there is $n_1=n_1(m)>0$ such that $ (\phi,\psi,\xi_0),
(\phi_n,\psi_n,\xi_{0n})\in E_{m+1}$ for all $n\ge n_1(m)$. Hence,
by Theorem~\ref{thLocalExistence}, there is $T\in(0,1]$ for which
problem~\eqref{eq1}, \eqref{eq3} has transverse topology
preserving solutions $(u,v)$ and $(u_n,v_n)$ with the
corresponding initial data. Moreover, any solution of
problem~\eqref{eq1}, \eqref{eq3} in $Q_{T}$ is transverse and
preserves topology.

We introduce the functions $a(t)$ and $a_n(t)$ corresponding to
$u$ and $u_n$ as described in part 3 of Lemma~\ref{ltNP0}. Then
the discontinuity points of the corresponding configuration
functions $\xi(x,t)$, $\xi_n(x,t)$ are given by
$b(t)=\max\limits_{s\in[0,t]}a(s)$ and
$b_n(t)=\max\limits_{s\in[0,t]}a_n(s)$.

2. Assume that there is $\varepsilon>0$ such that
\begin{equation}\label{eq1.1-4'}
\|u_n-u\|_{W_q^{2,1}(Q_{T};\cU)}\ge\varepsilon,\quad n=1,2,\dots,
\end{equation}
for some subsequence of $u_n$, which we denote $u_n$ again.
Theorem~\ref{thLocalExistence} implies that  $u_n$ and $a_n$ are
uniformly bounded in $W_q^{2,1}(Q_{T};\cU)$ and $C^\gamma[0,T]$,
respectively. Hence, we can choose subsequences of $u_n$ and $a_n$
(which we denote $u_n$ and $a_n$ again)  such that
\begin{gather}
\|u_n-\hat u\|_{C^\gamma(\oQ_{T};\cU)}\to 0,\quad \|(u_n)_x-\hat
u_x\|_{C^\gamma(\oQ_{T};\cU)}\to 0,\quad
n\to\infty,\label{eq1.1-8}
\\
\|a_n-\hat a\|_{C[0,T]}\to 0,\quad n\to\infty \label{eq1.1-9}
\end{gather}
for some function $\hat u\in C^\gamma(\oQ_{T};\cU)$ with  $\hat
u_x\in C^\gamma(\oQ_{T};\cU)$ and some  $\hat a\in C[0,T]$.

Set $\hat b(t)=\max\limits_{s\in[0,t]}\hat a(s)$. Due
to~\eqref{eq1.1-8}, \eqref{eq1.1-9}, and Lemma~\ref{lab}, we have
\begin{gather}
\label{eq1.1-11}
\|b_n-\hat b\|_{C[0,T]}\to 0,\quad n\to\infty,  \\
\label{eqContDepW}
W(\xi_0(x),\hat u(x,\cdot))(t)=\begin{cases} W_1(\hat u(x,t)), & 0\le x\le  \hat b(t),\\
W_2(\hat u(x,t)), & \hat b(t)<x\le 1.
\end{cases}
\end{gather}

3. Now we show that
\begin{equation}\label{eqContDepvn-v}
\sup\limits_{t\in[0,T]}\left(\|v_n(\cdot,t)-\hat
v(\cdot,t)\|_{L_q((0,1);\cV)}+\|v_{nt}(\cdot,t)-\hat
v_t(\cdot,t)\|_{L_q((0,1);\cV)}\right)\to 0,\quad n\to\infty,
\end{equation}
for some $\hat v$. Take an arbitrary $\delta>0$. It follows from
the assumptions of the theorem, from \eqref{eq1.1-8},
\eqref{eq1.1-11}, \eqref{eqContDepW}, and from Lemma~\ref{lu0b0v0}
that
\begin{equation}\label{eqContDepvn-v1}
\begin{gathered}
\|\psi_n-\psi_k\|_{L_q((0,1);\cV)}\le\delta,\quad
\|u_n(\cdot,t)-u_k(\cdot,t)\|_{L_q((0,1);\cU)}\le\delta,\\
\|W(\xi_n, u_n)(t)-W(\xi_k, u_k)(t)\|_{L_q((0,1);\cW)}\le\delta,
\end{gathered}
\end{equation}
provided $n,k$ are large enough. Estimates~\eqref{eqContDepvn-v1},
the second equation in~\eqref{eq1}, and the local Lipschitz
continuity of $g$ yield
$$
\|v_n(\cdot,t)-v_k(\cdot,t)\|_{L_q((0,1);\cV)}\le
(1+2L)\delta+L\int\limits_0^t
\|v_n(\cdot,s)-v_k(\cdot,s)\|_{L_q((0,1);\cV)}\,ds,
$$
where $L>0$ does not depend on $n,k$. Hence, by Gronwall's
inequality,
\begin{equation}\label{eqContDepvn-v2}
\|v_n(\cdot,t)-v_k(\cdot,t)\|_{L_q((0,1);\cV)}\le k_1\delta,
\end{equation}
where $k_1>0$ does not depend $\delta, n, k, $ and $t\in[0,T]$. A
similar inequality for the time derivative of $v_n$ follows
from~\eqref{eqContDepvn-v1}, \eqref{eqContDepvn-v2}, and from the
second equation in~\eqref{eq1}. Since $\delta>0$ is arbitrary,
\eqref{eqContDepvn-v} does hold.

Now consider~\eqref{eq1}, \eqref{eq3} with the subsequence
$\phi_n,\psi_n,\xi_{0n},u_n,v_n$ and pass   to the limit as
$n\to\infty$. Due to the uniqueness assumption, $(u,v)=(\hat
u,\hat v)$. Therefore,~\eqref{eq1.1-4'} is not true and we have
the  convergence for the whole sequence $(u_n,v_n)$.



\section{Uniqueness of solutions}\label{secUniqueness}

In this section, we prove Theorem \ref{tUniqueness}. For the
clarity of exposition, we restrict ourselves to the case where
initial data satisfy the equality $\gamma_{\alpha}(\varphi(\ob)) =
0$ in addition to Condition~\ref{condA1'}. (The case
$\gamma_{\alpha}(\varphi(\ob)) \ne 0$ can be treated easily
because then the hysteresis $W(\xi_0, u)$ remains constant on some
time interval.)

Set
$$
\ophi := \frac{1}{2}\frac{d}{dx}\gamma_{\alpha}(\phi(x))|_{x =
\ob} \ (> 0).
$$

We fix $T_1$ such that the conclusions of Lemma \ref{lOneRoot} are
true for the   $u$, $\hu$ on $(0, T_1)$. Let $a(t)$, $b(t)$, $\hat
a(t)$, $\hat b(t)$ be the functions defined in Lemma
\ref{lOneRoot} for     $u$ and $\hu$, respectively. We fix
$T\in(0,T_1)$ and $\delta>0$ such that the following hold for
$t\in[0,T]$$:$
\begin{gather}
\label{lem4.0.1}
 \frac{d}{dx}\gamma_{\alpha}(u(x,t))\ge\ophi, \quad x\in[\ob-\delta,\ob+\delta],\\
\label{lem4.0.5} \gamma_{\beta}(u(x,t)) < 0, \quad x\in[0,b(t)],
\end{gather}
and the analogous inequalities hold for $\hu$.

Due to~\eqref{lem4.0.1} and \eqref{lem4.0.5}, we have
\begin{equation}\label{eqHFreeBoundary}
\begin{aligned}
W(\xi_0(x),u(x,\cdot)(t)&=\begin{cases} W_1(u(x,t)), & 0\le x\le  b(t),\\
W_{-1}(u(x,t)), & b(t)<x\le 1,
\end{cases}\\
W(\xi_0(x),\hat u(x,\cdot)(t)&=\begin{cases} W_1(\hat u(x,t)), & 0\le x\le  \hat b(t),\\
W_{-1}(\hat u(x,t)), & \hat b(t)<x\le 1.
\end{cases}
\end{aligned}
\end{equation}

Let us now prove Theorem~\ref{tUniqueness}.

1. Denote $w=u-\hu$, $z = v - \hv$. The functions $w$, $z$ satisfy
the equations
\begin{equation}\label{eqLinearw}
\left\{\begin{aligned}
& w_t=w_{xx}+h_w(x,t),  \\
& z_t=h_z(x,t),
\end{aligned}\right.
\end{equation}
and the zero boundary and initial conditions, where
$$
\begin{aligned}
h_w(x,t)&= f(u, v, W(u)) - f(\hu, \hv, W(\hu)),\\
h_z(x,t)&= g(u, v, W(u)) - g(\hu, \hv, W(\hu)).
\end{aligned}
$$
Obviously, $h_{w, z}\in L_\infty(Q_T)$. The function $w$ can be
represented via the Green function $G(x,y,t,s)$ of the heat
equation with the Neumann boundary conditions:
\begin{equation*}
w(x,t)=\int\limits_0^t\int\limits_0^1 G(x,y,t,s)h_w(y,s)\,dy ds.
\end{equation*}

Therefore, using the estimate
$|G(x,y,t,s)|\le  \dfrac{k_1}{\sqrt{t-s}}$, $0<s<t$, with $k_1>0$
not depending on $(x,t)\in Q_T$ (see, e.g.,~\cite{Ivasisen}),
 we obtain
\begin{equation}\label{eqUniqueness1}
|w(x,t)|\le
k_1\int\limits_0^t\dfrac{ds}{\sqrt{t-s}}\int\limits_0^1
|h_w(y,s)|\,dy.
\end{equation}
Set $Z(t) = \int_0^1 |h_z(y, t)| dy$. Due to the second equation
in~\eqref{eqLinearw},
\begin{equation}\label{eqUniqueness1.1}
 Z(t) \leq \int\limits_0^t \int\limits_0^1  |h_z(y, s)| dy ds.
\end{equation}

2. Now we prove that, for some   $k_2 > 0$,
\begin{equation}\label{eqUniqueness1.2}
\int\limits_0^1 |h_{w, z}(y, s)| dy \leq k_2(\|w\|_{C(\oQ_T)} +
\|Z\|_{L_\infty(0, T)}), \quad s \in (0, T).
\end{equation}

Let us prove this inequality for the function $h_w$, assuming that
$b(s)<\hat b(s)$. (The cases of $h_z$ and $b(s)\ge\hat b(s)$ are
treated analogously.) Since $f$ is locally Lipschitz,
\begin{multline}\label{eqUniqueness1.3}
\int\limits_0^1 |h_w(y, s)| dy \leq k_3 \int\limits_0^1 (|w(y, s)| + |z(y, s)| + |W(u(y, s)) - W(\hu(y, s))|) \leq \\
\leq k_3\left(\|w\|_{C(\oQ_T)} + \|Z\|_{L_\infty(0, T)} +
\int\limits_0^1 |W(u(y, s)) - W(\hu(y, s))| dy\right),
\end{multline}
where $k_3>0$ and the constants $k_4, k_5, \ldots >0$ below do not
depend on $s\in[0,T]$.

Denote $\theta(y, s) = W(u(y, s)) - W(\hu(y, s))$. Due
to~\eqref{eqHFreeBoundary}, we have
$$
\theta(y, s)=
\begin{cases}
W_1(u)-W_1(\hat u),& 0<y< b(s),\\
W_{-1}(u)-W_1(\hat u),& b(s)<y< \hat b(s),\\
W_{-1}(u)-W_{-1}(\hat u),& \hat b(s)<y< 1.
\end{cases}
$$

2.1. Inequality \eqref{lem4.0.5} implies that
$\gamma_{\beta}(u(y,s))<0$,$\gamma_{\beta}(\hat u(y,s))<0$ on the
closed set $\{(y,s): y\in[0, b(s)],\ s\in[0,T]\}$. Hence, the
values $\gamma_{\beta}(u(y,s))$ and $\gamma_{\beta}(\hat u(y,s))$
are separated from 0. Therefore, using Condition~\ref{condUniq},
we obtain
\begin{equation}\label{eqUniqueness2}
\int\limits_0^{b(s)} |\theta(y,s)|\,dy \le k_4
\int\limits_0^{b(s)}|u(y,s)-\hat u(y,s)|\,dy \le
k_4\|w\|_{C(\oQ_T)}.
\end{equation}

2.2. Boundedness of $W_1(\hat u)$ and $W_{-1}(u)$ for
$(y,s)\in\oQ_T$ and Lemma \ref{lab} imply
\begin{equation}\label{eqUniqueness3'}
\int\limits_{b(s)}^{\hat b(s)} |\theta(y,s)|\,dy\le k_5
\int\limits_{b(s)}^{\hat b(s)}\,dy\le k_5\|b-\hat b\|_{C[0,T]}\le
k_5 \|a-\hat a\|_{C[0,T]}.
\end{equation}
 Using \eqref{lem4.0.1}, we obtain for any
$t\in[0,T]$ the inequalities
\begin{equation}\label{eqUniqueness4'}
\begin{aligned}
|a(t)-\hat a(t)| &\le
\dfrac{1}{\ophi}|\gamma_{\alpha}(u(a(t),t))-\gamma_{\alpha}(\hat
u(a(t),t))| \le\\
&\le \dfrac{L_{\alpha}}{\ophi}|u(a(t),t)-\hat u(a(t),t)| \le
\dfrac{L_{\alpha}}{\ophi}\|u-\hat u\|_{C(\oQ_T)},
\end{aligned}
\end{equation}
where $L_{\alpha} > 0$ is a respective Lipschitz constant for
$\gamma_\alpha(u)$ and hence does not depend on $T \in (0, T_1)$.
Inequalities~\eqref{eqUniqueness3'} and~\eqref{eqUniqueness4'}
yield
\begin{equation}\label{eqUniqueness3}
\int\limits_{b(s)}^{\hat b(s)} |\theta(y,s)|\,dy\le k_6
\|w\|_{C(\oQ_T)}.
\end{equation}

2.3. Let $y \in [\hat b(s), \ob+\delta]$. Inequality
\eqref{lem4.0.1} and the mean-value theorem imply
$$
\begin{gathered}
\gamma_{\alpha}(\hu(y, s)) = \gamma_{\alpha}(\hu(y, s)) -
\gamma_{\alpha}(\hu(\hat{a}(s), s)) \ge (y-\hat
a(s))\ophi\ge (y-\hat b(s))\ophi,\\
|\gamma_{\alpha}(u(y,s))| \ge (y-b(s))\ophi.
\end{gathered}
$$
Taking into account these two inequalities and using
Condition~\ref{condUniq}, we obtain
\begin{equation}\label{eqUniqueness4}
\int\limits_{\hat b(s)}^{\ob+\delta} |\theta(y,s)|\,dy \le k_7
\int\limits_{\hat b(s)}^{\ob+\delta} \dfrac{|u(y,s)-\hat
u(y,s)|}{(y-\hat b(s))^\sigma}\,dy \le k_8\|w\|_{C(\oQ_T)}.
\end{equation}

2.4. Similarly to item 2.1, we conclude that
\begin{equation}\label{eqUniqueness5}
\int\limits_{\ob+\delta}^1 |\theta(y,s)|\,dy \le
k_9\|w\|_{C(\oQ_T)}.
\end{equation}

Finally, \eqref{eqUniqueness1.3}--\eqref{eqUniqueness5} imply
\eqref{eqUniqueness1.2}.

3. Combining
estimates~\eqref{eqUniqueness1}--\eqref{eqUniqueness1.2}, we
obtain
$$
\begin{aligned}
|w(x,t)| & \le k_{10} (\|w\|_{C(\oQ_T)} + \|Z\|_{L_\infty(0, T)})
\int\limits_0^t \dfrac{ds}{\sqrt{t-s}}\\
&=2 k_{10} T^{1/2}(\|w\|_{C(\oQ_T)} + \|Z\|_{L_\infty(0, T)}),
\\
Z(t) &\le k_2 T(\|w\|_{C(\oQ_T)} + \|Z\|_{L_\infty(0, T)}).
\end{aligned}
$$
Taking the supremum with respect to $t\in(0,T)$, we see that
$$
\|w\|_{C(\oQ_T)} + \|Z\|_{L_\infty(0, T)} \leq (2 k_{10} T^{1/2} +
k_2 T)(\|w\|_{C(\oQ_T)} + \|Z\|_{L_\infty(0, T)}).
$$
Thus, $w=0$ and $z = 0$, provided that $T>0$ is small enough. \qed

{\small {\em Authors' addresses}: {\em Pavel Gurevich}, Free
University Berlin, Arnimallee 3, Berlin, 14195, Germany; Peoples'
Friendship University, Mikluho-Maklaya str. 6, Moscow, 117198,
Russia; e-mail: \texttt{gurevichp@\allowbreak gmail.com}. {\em
Sergey Tikhomirov}, Chebyshev Laboratory, Saint-Petersburg State
University, 14th line of Vasilievsky island, 29B,
Saint-Petersburg, 199178, Russia; Max Planck Institute for
Mathematics in the Sciences, Inselstrasse 22, 04103, Leipzig,
Germany; e-mail: \texttt{sergey.tikhomirov@\allowbreak gmail.com}.
}


{\small
\begin{thebibliography}{999}


\bibitem{AikiKopfova}
{\it T. Aiki, J. Kopfov\'a}. A Mathematical Model for Bacterial
Growth Described by a Hysteresis Operator. Recent advances in
nonlinear analysis -- Proceedings of the International Conference
on Nonlinear Analysis. Hsinchu, Taiwan, 2006.

\bibitem{Alt}
{\it H. W. Alt}. On the thermostat problem. Control Cyb., {\bf 14}, 171--193 (1985).

\bibitem{UrAp}
\textit{D. Apushinskaya, H.Shahgholian, N.N.Uraltseva.} On the Lipschitz property of the free boundary in a parabolic problem with an obstacle.
St. Petersburg Math. J., \textbf{15}, no. 3, 375--391  (2004).

\bibitem{GurTikhSIAM13}
{\it P. Gurevich, R. Shamin, S. Tikhomirov.} Reaction-Diffusion Equations with Spatially Distributed Hysteresis. SIAM J. Math. Anal. {\bf 45}, no. 3, 1328-1355 (2013).

\bibitem{GurTikhNonlinAnal}
{\it P. Gurevich, S. Tikhomirov.} Uniqueness of transverse solutions for reaction-diffusion equations with spatially distributed hysteresis. Nonlinear Analysis Series A:
Theory, Methods and Applications, \textbf{75}, 6610--6619 (2012).

\bibitem{Jaeger1}
{\it F. C. Hoppensteadt, W. J\"ager}.  Lecture Notes in
Biomathematics {\bf 38}, 68--81 (1980).

\bibitem{Jaeger2}
{\it F. C. Hoppensteadt, W. J\"ager, C. Poppe}. Modelling of
Patterns in Space and Time, Lecture Notes in Biomath. {\bf 55}
(Springer), 123--134 (1984).

\bibitem{Ilin}
{\it A. M. Il'in, B. A. Markov}. A nonlinear diffusion equation and Liesegang rings. Doklady Mathematics, {\bf 84}, No. 2, 730--733 (2011).

\bibitem{Ivasisen} {\it S. D. Ivasi$\check{s}$en}  Green's matrices of boundary value problems for Petrovski parabolic systems of general form, II,
Mat. Sb., {\bf 114(156)}, No. 4, 523--565 (1981); English transl.
in Math. USSR Sbornik, {\bf 42} 461--489 (1982).

\bibitem{Klein}
{\it O. Klein}.   Representation of hysteresis operators acting on
vector-valued monotaffine functions, Adv. Math. Sci. Appl., {\bf
22}, 471--500 (2012).

\bibitem{Kopfova}
{\it J. Kopfov\'a}. Hysteresis in biological models. Proceedings
of the conference ``International Workshop on Multi-rate
processess and hysteresis'', Journal of Physics, Conference
Series, {\bf 55}, 130--134 (2006).

\bibitem{KrasnBook}
{\it M. A. Krasnosel'skii, A. V. Pokrovskii}.  Systems with
Hysteresis. Springer-Verlag. Berlin--Heidelberg--New York (1989).

\bibitem{LadSolUral} {\it O. A. Ladyzhenskaya, V. A. Solonnikov, N. N.
Uraltseva}.  Linear and Quasilinear Equations of Parabolic Type.
Nauka, Moscow, 1967.

\bibitem{Rothe}
{\it F. Rothe}.  Global Solutions of Reaction-Diffusion Systems.
Lecture Notes in Mathematics, 1072. Springer-Verlag, Berlin
(1984).

\bibitem{ShUr}
\textit{H. Shahgholian, N. Uraltseva, G. S. Weiss.} A parabolic two-phase obstacle-like equation. Advances in
Mathematics, {\bf 221}, 861--881  (2009).

\bibitem{Smoller}
{\it J. Smoller}.  Shock Waves and Reaction-Diffusion Equations.
Second edition. Grundlehren der Mathematischen Wissenschaften
[Fundamental Principles of Mathematical Sciences], 258.
Springer-Verlag, New York (1994).


\bibitem{VisintinSpatHyst}
{\it A. Visintin}. Evolution problems with hysteresis in the source term. SIAM J. Math. Anal, {\bf 17}, 1113--1138
(1986).

\bibitem{Visintin}
{\it A. Visintin}.   Differential Models of Hysteresis.
Springer-Verlag. Berlin --- Heidelberg (1994).

\end{thebibliography}}
\end{document}